# 32
# Biology-Derived Algorithms in Engineering Optimization

Xin-She Yang

## 32.1 Introduction

Biology-derived algorithms are an important part of computational sciences, which are essential to many scientific disciplines and engineering applications. Many computational methods are derived from or based on the analogy to natural evolution and biological activities, and these biologically inspired computations include genetic algorithms, neural networks, cellular automata, and other algorithms. However, a substantial amount of computations today are still using conventional methods such as finite difference, finite element, and finite volume methods. New algorithms are often developed in the form of a hybrid combination of biology-derived algorithms and conventional methods, and this is especially true in the field of engineering optimizations. Engineering problems with optimization objectives are often difficult and time consuming, and the application of nature or biology-inspired algorithms in combination with the conventional optimization methods has been very successful in the last several decades.

There are five paradigms of nature-inspired evolutionary computations: genetic algorithms, evolutionary programming, evolutionary strategies, genetic programming, and classifier systems (Holland, 1975; Goldberg, 1989; Mitchell, 1996; Flake, 1998). Genetic algorithm (GA), developed by John Holland and his collaborators in the 1960s and 1970s, is a model or abstraction of biological evolution, which includes the following operators: crossover, mutation, inversion, and selection. This is done by the representation within a computer of a population of individuals corresponding to chromosomes in terms of a set of character strings, and the individuals in the population then evolve through the crossover and mutation

**32**-585

of the string from parents, and the selection or survival according to their fitness. Evolutionary programming (EP), first developed by Lawrence J. Fogel in 1960, is a stochastic optimization strategy similar to GAs. But it differs from GAs in that there is no constraint on the representation of solutions in EP and the representation often follows the problem. In addition, the EPs do not attempt to model genetic operations closely in the sense that the crossover operation is not used in EPs. The mutation operation simply changes aspects of the solution according to a statistical distribution, such as multivariate Gaussian perturbations, instead of bit-flopping, which is often done in GAs. As the global optimum is approached, the rate of mutation is often reduced. Evolutionary strategies (ESs) were conceived by Ingo Rechenberg and Hans-Paul Schwefel in 1963, later joined by Peter Bienert, to solve technical optimization problems. Although they were developed independently of one another, both ESs and EPs have many similarities in implementations. Typically, they both operate on real-values to solve real-valued function optimization in contrast with the encoding in Gas. Multivariate Gaussian mutation with zero mean are used for each parent population and appropriate selection criteria are used to determine which solution to keep or remove. However, EPs often use stochastic selection via a tournament and the selection eliminates those solutions with the least wins, while the ESs use deterministic selection criterion that removes the worst solutions directly based on the evaluations of certain functions (Heitkotter and Beasley, 2000). In addition, recombination is possible in an ES as it is an abstraction of evolution at the level of individual behavior in contrast to the abstraction of evolution at the level of reproductive populations and no recombination mechanisms in EPs.

The aforementioned three areas have the most impact in the development of evolutionary computations, and, in fact, evolutionary computation has been chosen as the general term that encompasses all these areas and some new areas. In recent years, two more paradigms in evolutionary computation have attracted substantial attention: Genetic programming and classifier systems. Genetic programming (GP) was introduced in the early 1990s by John Koza (1992), and it extends GAs using parse trees to represent functions and programs. The programs in the population consist of elements from the function sets, rather than fixed-length character strings, selected appropriately to be the solutions to the problems. The crossover operation is done through randomly selected subtrees in the individuals according to their fitness; the mutation operator is not used in GP. On the other hand, a classifier system (CFS), another invention by John Holland, is an adaptive system that combines many methods of adaptation with learning and evolution. Such hybrid systems can adapt behaviors toward a changing environment by using GAs with adding capacities such as memory, recursion, or iterations. In fact, we can essentially consider the CFSs as general-purpose computing machines that are modified by both environmental feedback and the underlying GAs (Holland, 1975, 1996; Flake, 1998).

AQ: Please check if the insertion of (1992) here is ok.

AQ: Holland, 1996, is not listed in the References. Is it 1995 as in the list?

Biology-derived algorithms are applicable to a wide variety of optimization problems. For example, optimization functions can have discrete, continuous, or even mixed parameters without any a priori assumptions about their continuity and differentiability. Thus, evolutionary algorithms are particularly suitable for parameter search and optimization problems. In addition, they are easy for parallel implementation. However, evolutionary algorithms are usually computationally intensive, and there is no absolute guarantee for the quality of the global optimizations. Besides, the tuning of the parameters can be very difficult for any given algorithms. Furthermore, there are many evolutionary algorithms with different suitabilities and the best choice of a particular algorithm depends on the type and characteristics of the problems concerned. However, great progress has been made in the last several decades in the application of evolutionary algorithms in engineering optimizations. In this chapter, we will focus on some of the important areas of the application of GAs in engineering optimizations.

## 32.2 Biology-Derived Algorithms

There are many biology-derived algorithms that are popular in evolutionary computations. For engineering applications in particular, four types of algorithms are very useful and hence relevant. They are GAs, photosynthetic algorithms (PAs), neural networks, and cellular automata. We will briefly discuss these algorithms in this section, and we will focus on the application of GAs and PAs in engineering optimizations in Section 32.3.

### 32.2.1 Genetic Algorithms

The essence of GAs involves the encoding of an optimization function as arrays of bits or character strings to represent the chromosomes, the manipulation operations of strings by genetic operators, and the selection according to their fitness to find a solution to the problem concerned. This is often done by the following procedure: (1) encoding of the objectives or optimization functions; (2) defining a fitness function or selection criterion; (3) creating a population of individuals; (4) evolution cycle or iterations by evaluating the fitness of all the individuals in the population, creating a new population by performing crossover, mutation, and inversion, fitness-proportionate reproduction, etc., and replacing the old population and iterating using the new population; (5) decoding the results to obtain the solution to the problem.

One iteration of creating new populations is called a generation. Fixed-length character strings are used in most GAs during each generation although there is substantial research on variable-length string and coding structures. The coding of objective functions is usually in the form of binary arrays or real-valued arrays in adaptive GAs. For simplicity, we use the binary string for coding for describing genetic operators. Genetic operators include crossover, mutation, inversion, and selection. The crossover of two parent strings is the main operator with highest probability $p_c$ (usually, 0.6 to 1.0) and is carried out by switching one segment of one string with the corresponding segment on another string at a random position (see Figure 32.1). The crossover carried out in this way is a single-point crossover. Crossover at multiple points is also used in many GAs to increase the efficiency of the algorithms. The mutation operation is achieved by the flopping of randomly selected bits, and the mutation probability $p_m$ is usually small (say, 0.001 to 0.05), while the inversion of some part of a string is done by interchanging 0 and 1. The selection of an individual in a population is carried out by the evaluation of its fitness, and it can remain in the new generation if a certain threshold of fitness is reached or the reproduction of a population is fitness-proportionate. One of the key parts is the formulation or choice of fitness functions that determines the selection criterion in a particular problem.

Further, just as crossover can be carried out at multiple points, mutations can also occur at multiple sites. More complex and adaptive GAs are being actively researched and there is a vast literature on this topic. In Section 32.3, we will give examples of GAs and their applications in engineering optimizations.

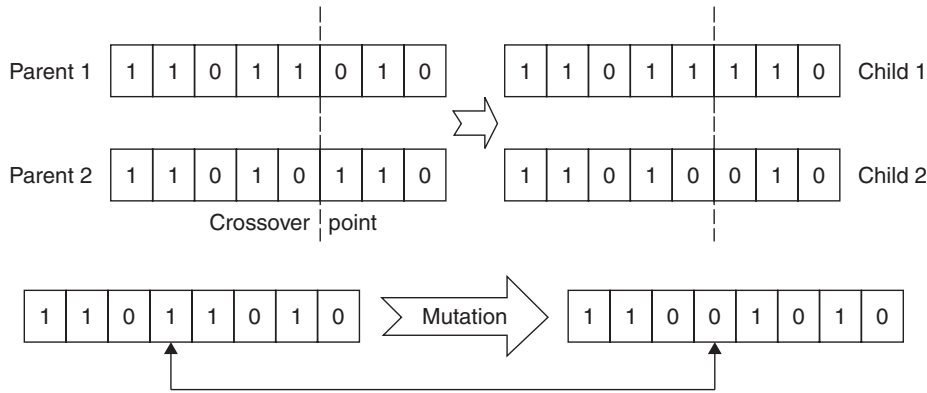

**FIGURE 32.1**   Diagram of crossover and mutation in a GA.

### 32.2.2 Photosynthetic Algorithms

The PA was first introduced by Murase (2000) to optimize the parameter estimation in finite element inverse analysis. The PA is a good example of biology-derived algorithms in the sense that its computational procedure corresponds well to the real photosynthesis process in green plants. Photosynthesis uses water and $CO_2$ to produce glucose and oxygen when there is light and in the presence of chloroplasts. The overall reaction

> AQ: Murose is changed to Murase as per the Reference list. Trust this is OK.

$$6CO_2 + 12H_2O \xrightarrow{\text{light and green plants}} C_6H_{12}O_6 + 6O_2 + 6H_2O$$

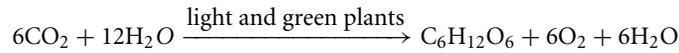

is just a simple version of a complicated process. Other factors, such as temperature, concentration of $CO_2$, water content, etc., being equal, the reaction efficiency depends largely on light intensity. The important part of photosynthetic reactions is the dark reactions that consist of a biological process including two cycles: the Benson–Calvin cycle and photorespiration cycle. The balance between these two cycles can be considered as a natural optimization procedure that maximizes the efficiency of sugar production under the continuous variations of light energy input (Murase, 2000).

Murase's PA uses the rules governing the conversion of carbon molecules in the Benson–Calvin cycle (with a product or feedback from DHAP) and photorespiration reactions. The product DHAP serves as the knowledge strings of the algorithm and optimization is reached when the quality or the fitness of the products no longer improves. An interesting feature of such algorithms is that the stimulation is a function of light intensity that is randomly changing and affects the rate of photorespiration. The ratio of $O_2$ to $CO_2$ concentration determines the ratio of the Benson–Calvin and photorespiration cycles. A PA consists of the following steps: (1) the coding of optimization functions in terms of fixed-length DHAP strings (16-bit in Murase's PA) and a random generation of light intensity ($L$); (2) the $CO_2$ fixation rate $r$ is then evaluated by the following equation: $r = V_{max}/(1 + A/L)$, where $V_{max}$ is the maximum fixation rate of $CO_2$ and $A$ is its affinity constant; (3) either the Benson–Calvin cycle or photorespiration cycle is chosen for the next step, depending on the $CO_2$ fixation rate, and the 16-bit strings are shuffled in both cycles according to the rule of carbon molecule combination in photosynthetic pathways; (4) after some iterations, the fitness of the intermediate strings is evaluated, and the best fit remains as a DHAP, then the results are decoded into the solution of the optimization problem (see Figure 32.2).

> AQ: Please provide expansion for DHAP.

In the next section, we will present an example of parametric inversion and optimization using PA in finite element inverse analysis.

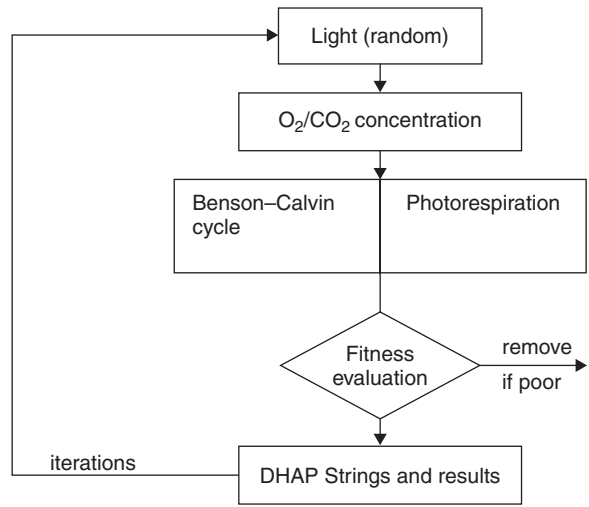

**FIGURE 32.2** Scheme of Murase's PAs.

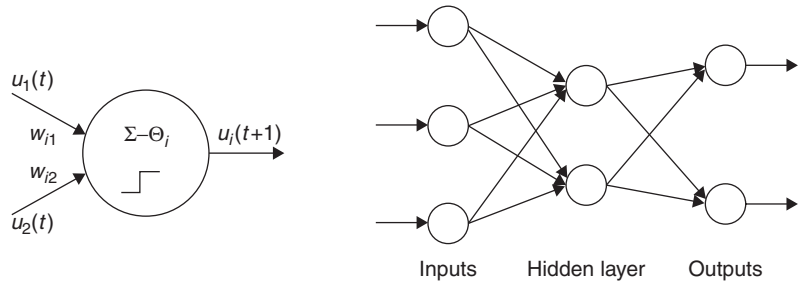

**FIGURE 32.3** Diagram of a McCulloch–Pitts neuron (left) and neural networks (right).

### 32.2.3 Neural Networks

Neural networks, and the associated machine-learning algorithms, is one more type of biology-inspired algorithms, which uses a network of interconnected neurons with the intention of imitating the neural activities in human brains. These neurons have high interconnectivity and feedback, and their connectivity can be weakened or enhanced during the learning process. The simplest model for neurons is the McCulloch–Pitts model (see Figure 32.3) for multiple inputs $u_1(t), \ldots, u_n(t)$, and the output $u_i(t)$. The activation in a neuron or a neuron's state is determined by

$$u_i(t+1) = H\left[\sum_{j=1}^{n} w_{ij} u_i(t) - \Theta_i\right],$$

where $H(x)$ is the Heaviside unit step function that $H(x) = 1$ if $x \geq 0$ and $H(x) = 0$ otherwise. The weight coefficient $w_{ij}$ is considered as the synaptic strength of the connection from neuron $j$ to neuron $i$. For each neuron, it can be activated only if its threshold $\Theta_i$ is reached. One can consider a single neuron as a simple computer that gives the output 1 or yes if the weighted sum of incoming signals is greater than the threshold, otherwise it outputs 0 or no.

Real power comes from the combination of nonlinear activation functions with multiple neurons (McCulloch and Pitts, 1943; Flake, 1998). Figure 32.3 also shows an example of feed-forward neural networks. The key element of an artificial neural network (ANN) is the novel structure of such an information processing system that consists of a large number of interconnected processing neurons. These neurons work together to solve specific problems by adaptive learning through examples and self-organization. A trained neural network in a given category can solve and answer what-if type questions to a particular problem when new situations of interest are given. Due to the real-time capability and parallel architecture as well as adaptive learning, neural networks have been applied to solve many real-world problems such as industrial process control, data validation, pattern recognition, and other systems of artificial intelligence such as drug design and diagnosis of cardiovascular conditions. Optimization is just one possibility of such applications, and often the optimization functions can change with time as is the case in industrial process control, target marketing, and business forecasting (Haykins, 1994). On the other hand, the training of a network may take considerable time and a good training database or examples that are specific to a particular problem are required. However, neural networks will, gradually, come to play an important role in engineering applications because of its flexibility and adaptability in learning.

### 32.2.4 Cellular Automata

Cellular automata (CA) were also inspired by biological evolution. On a regular grid of cells, each cell has finite number of states. Their states are updated according to certain local rules that are functions of the states of neighbor cells and the current state of the cell concerned. The states of the cells evolve with time in a discrete manner, and complex characteristics can be observed and studied. For more details on this topic, readers can refer to Chapters 1 and 22 in this handbook.

There is some similarity between finite state CA and conventional numerical methods such as finite difference methods. If one considers the finite different method as real-valued CA, and the real-values are always converted to finite discrete values due to the round-off in the implementation on a computer, then there is no substantial difference between a finite difference method and a finite state CA.

However, CA are easier to parallelize and more numerically stable. In addition, finite difference schemes are based on differential equations and it is sometimes straightforward to formulate a CA from the corresponding partial differential equations via appropriate finite differencing procedure; however, it is usually very difficult to conversely obtain a differential equation for a given CA (see Chapter 22 in this handbook).

An optimization problem can be solved using CA if the objective functions can be coded to be associated with the states of the CA and the parameters are properly associated with automaton rules. This is an area under active research. One of the advantages of CA is that it can simulate many processes such as reaction–diffusion, fluid flow, phase transition, percolation, waves, and biological evolution. Artificial intelligence also uses CA intensively.

### 32.2.5 Optimization

Many problems in engineering and other disciplines involve optimizations that depend on a number of parameters, and the choice of these parameters affects the performance or objectives of the system concerned. The optimization target is often measured in terms of objective or fitness functions in qualitative models. Engineering design and testing often require an iteration process with parameter adjustment. Optimization functions are generally formulated as:

$$\text{Optimize: } f(\mathbf{x}),$$
$$\text{Subject to: } g_i(\mathbf{x}) \geq 0, \quad i = 1, 2, \ldots, N; \quad h_j(x) = 0, \quad j = 1, 2, \ldots, M.$$
$$\text{where } \mathbf{x} = (x_1, x_2, \ldots, x_n), \quad \mathbf{x} \in \Omega(\text{parameter space}).$$

Optimization can be expressed either as maximization or more often as minimization (Deb, 2000). As parameter variations are usually very large, systematic adaptive searching or optimization procedures are required. In the past several decades, researchers have developed many optimization algorithms. Examples of conventional methods are hill climbing, gradient methods, random walk, simulated annealing, heuristic methods, etc. Examples of evolutionary or biology-inspired algorithms are GAs, photosynthetic methods, neural network, and many others.

The methods used to solve a particular problem depend largely on the type and characteristics of the optimization problem itself. There is no universal method that works for all problems, and there is generally no guarantee to find the optimal solution in global optimizations. In general, we can emphasize on the best estimate or suboptimal solutions under the given conditions. Knowledge about the particular problem concerned is always helpful to make the appropriate choice of the best or most efficient methods for the optimization procedure. In this chapter, however, we focus mainly on biology-inspired algorithms and their applications in engineering optimizations.

## 32.3 Engineering Optimization and Applications

Biology-derived algorithms such as GAs and PAs have many applications in engineering optimizations. However, as we mentioned earlier, the choice of methods for optimization in a particular problem depends on the nature of the problem and the quality of solutions concerned. We will now discuss optimization problems and related issues in various engineering applications.

### 32.3.1 Function and Multilevel Optimizations

For optimization of a function using GAs, one way is to use the simplest GA with a fitness function: $F = A - y$ with $A$ being the large constant and $y = f(\mathbf{x})$, thus the objective is to maximize the fitness function and subsequently minimize the objective function $f(\mathbf{x})$. However, there are many different ways of defining a fitness function. For example, we can use the individual fitness assignment relative to the whole population

$$F(x_i) = \frac{f(x_i)}{\sum_{i=1}^{N} f(x_i)},$$

where $x_i$ is the phenotypic value of individual $i$, and $N$ is the population size. For the generalized De Jong's (1975) test function

$$f(\mathbf{x}) = \sum_{i=1}^{n} x^{2\alpha}, \quad |x| \leq r, \quad \alpha = 1, 2, \ldots, m,$$

where $\alpha$ is a positive integer and $r$ is the half-length of the domain. This function has a minimum of $f(\mathbf{x}) = 0$ at $\mathbf{x} = 0$. For the values of $\alpha = 3, r = 256$, and $n = 40$, the results of optimization of this test function are shown in Figure 32.4 using GAs.

AQ: Please check the change of $f(x)$ to $f(\mathbf{x})$

The function we just discussed is relatively simple in the sense that it is single-peaked. In reality, many functions are multi-peaked and the optimization is thus multileveled. Keane (1995) studied the following bumby function in a multi-peaked and multileveled optimization problem

$$f(x, y) = \frac{\sin^2(x - y) \sin^2(x + y)}{\sqrt{x^2 + y^2}}, \quad 0 < x, \ y < 10.$$

The optimization problem is to find $(x, y)$ starting $(5, 5)$ to maximize the function $f(x, y)$ subject to: $x + y \leq 15$ and $xy \geq 3/4$. In this problem, optimization is difficult because it is nearly symmetrical about $x = y$, and while the peaks occur in pairs one is bigger than the other. In addition, the true maximum is $f(1.593, 0.471) = 0.365$, which is defined by a constraint boundary. Figure 32.5 shows the surface variation of the multi-peaked bumpy function.

Although the properties of this bumpy function make it difficult for most optimizers and algorithms, GAs and other evolutionary algorithms perform well for this function and it has been widely used as a test function in GAs for comparative studies of various evolutionary algorithms or in multilevel optimization environments (El-Beltagy and Keane, 1999).

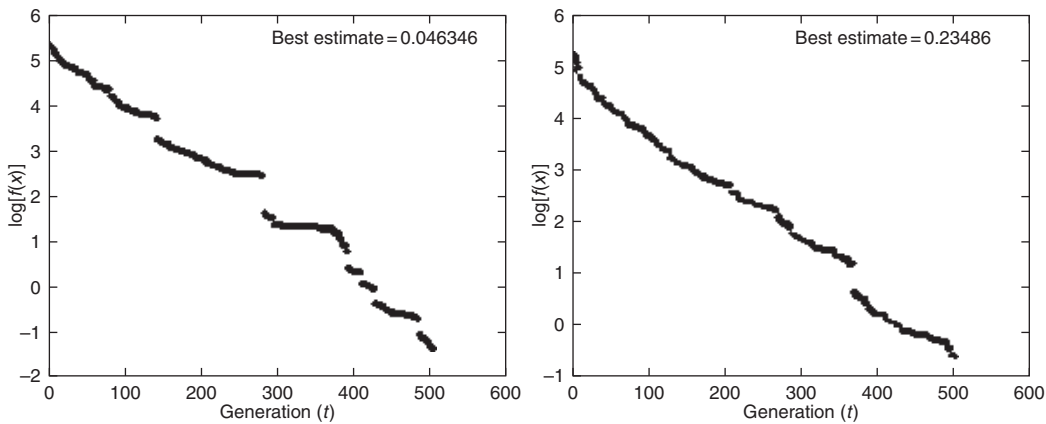

**FIGURE 32.4** Function optimization using GAs. Two runs will give slightly different results due to the stochastic nature of GAs, but they produce better estimates: $f(\mathbf{x}) \to 0$ as the generation increases.

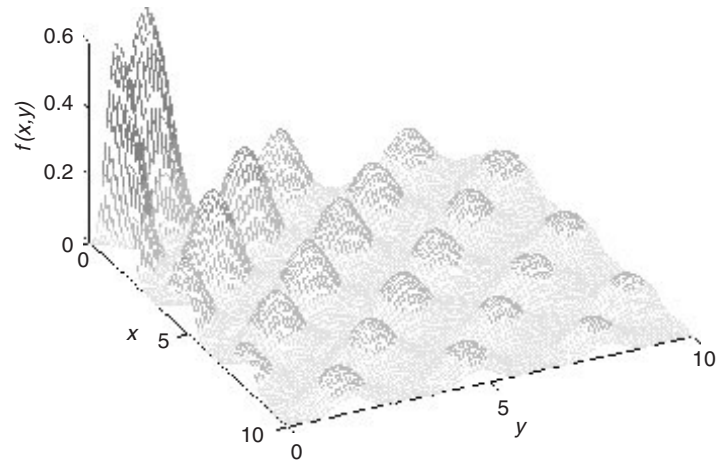

**FIGURE 32.5**　Surface of the multi-peaked bumpy function.

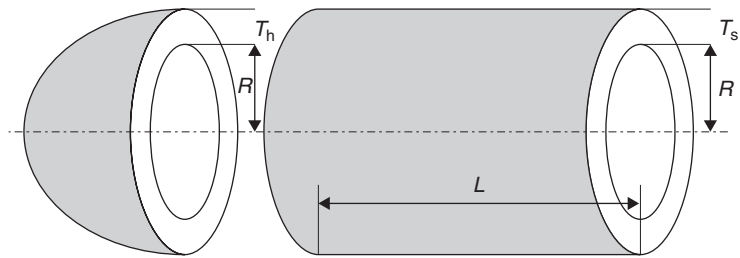

**FIGURE 32.6**　Diagram of the pressure vessel.

## 32.3.2　Shape Design and Optimization

Most engineering design problems, especially in shape design, aim to reduce the cost, weight, and volume and increase the performance and quality of the products. The optimization process starts with the transformation of design specification and descriptions into optimization functions and constraints. The structure and parameters of a product depend on the functionality and manufacturability, and thus considerable effort has been put into the modeling of the design process and search technique to find the optimal solution in the search space, which comprises the set of all designs with all allowable values of design parameters (Renner and Ekart, 2003). Genetic algorithms have been applied in many areas of engineering design such as conceptual design, shape optimization, data fitting, and robot path design.

A well-studied example is the design of a pressure vessel (Kannan and Kramer, 1994; Coello, 2000) using different algorithms such as augmented Lagrangian multiplier and GAs. Figure 32.6 shows the diagram of the parameter notations of the pressure vessel. The vessel is cylindrical and capped at both ends by hemispherical heads with four design variables: thickness of the shell $T_s$, thickness of the head $T_h$, inner radius $R$, and the length of the cylindrical part $L$. The objective of the design is to minimize the total cost including that of the material, forming, as well as welding. Using the notation given by Kannan and

Kramer, the optimization problem can be expressed as:

$$\text{Minimize: } f(\mathbf{x}) = 0.6224 x_1 x_2 x_3 + 1.7781 x_2 x_3^2 + 3.1611 x_1^2 x_4 + 19.84 x_1^2 x_3.$$

$$\mathbf{x} = (x_1, x_2, x_3, x_4)^{\text{T}} = (T_s, T_h, R, L)^{\text{T}},$$

$$\text{Subject to: } g_1(\mathbf{x}) = -x_1 + 0.0193 x_3 \leq 0, \quad g_2(\mathbf{x}) = -x_2 + 0.00954 x_3 \leq 0,$$

$$g_3(\mathbf{x}) = -\pi x_3^2 x_4 - 4\pi x_3^3/3 + 1296000 \leq 0, \quad g_4(\mathbf{x}) = x_4 - 240 \leq 0.$$

The values for $x_1, x_2$ should be considered as integer multiples of 0.0625. Using the same constraints as given in Coello (2000), the variables are in the ranges: $1 \leq x_1, x_2 \leq 99, 10.0000 \leq x_3. x_4 \leq 100.0000$ (with a four-decimal precision). By coding the GAs with a population of 44-bit strings for each individual (4-bits for $x_1 \cdot x_2$; 18-bits for $x_3 \cdot x_4$), similar to that by Wu and Chow (1994), we can solve the optimization problem for the pressure vessel. After several runs, the best solution obtained is $x_* = (1.125, 0.625, 58.2906, 43.6926)$ with $f(\mathbf{x}) = 7197.9912\$$, which is compared with the results $x_* = (1.125, 0.625, 28.291, 43.690)$ and $f(\mathbf{x}) = 7198.0428\$$ obtained by Kannan and Kramer (1994).

### 32.3.3 Finite Element Inverse Analysis

The usage and efficiency of Murase's PA described in Section 32.3.4 can be demonstrated in the application of finite element inverse analysis. Finite element analysis (FEA) in structural engineering is forward modeling as the aims are to calculate the displacements at various positions for given loading conditions and material properties such as Young's modulus ($E$) and Poisson's ratio ($\nu$). This forward FEA is widely used in engineering design and applications. Sometimes, inverse problems need to be solved. For a given structure with known loading conditions and measured displacement, the objective is to find or invert the material properties $E$ and $\nu$, which may be required for testing new materials and design optimizations. It is well known that inverse problems are usually very difficult, and this is especially true for finite element inverse analysis.

To show how it works, we use a test example similar to that proposed by Murase (2000). A simple beam system of 5 unit length × 10 unit length (see Figure 32.7) consists of five nodes and four elements whose Young's modulus and Poisson's ratio may be different. Nodes 1 and 2 are fixed, and a unit vertical load is applied at node 4 while the other nodes deform freely. Let us denote the displacement vector

$$\mathbf{U} = (u_1, v_1, u_2, v_2, u_3, v_3, u_4, v_4, u_5, v_5)^{\text{T}},$$

where $u_1 = v_1 = u_2 = v_2 = 0$ (fixed). Measurements are made for other displacements. By using the PA with the values of the $CO_2$ affinity $A = 10000$, light intensity $L = 10^4$ to $5 \times 10^4$ lx, and maximum $CO_2$ fixation speed $V_{\max} = 30$, each of eight elastic modulus $(E_i, v_i)(i = 1, 2, 3, 4)$ is coded as a 16-bit DHAP molecule string. For a target vector

$$\mathbf{Y} = (E_1, v_1, E_2, v_2, E_3, v_3, E_4, v_4) = (600, 0.25, 400, 0.35, 450, 0.30, 350, 0.32),$$

and measured displacements $\mathbf{U} = (0, 0, 0, 0, -0.0066, -0.0246, 0.0828, -0.2606, 0.0002, -0.0110)$, the best estimates after 500 iterations from the optimization by the PA are $\mathbf{Y} = (580, 0.24, 400, 0.31, 460, 0.29, 346, 0.26)$.

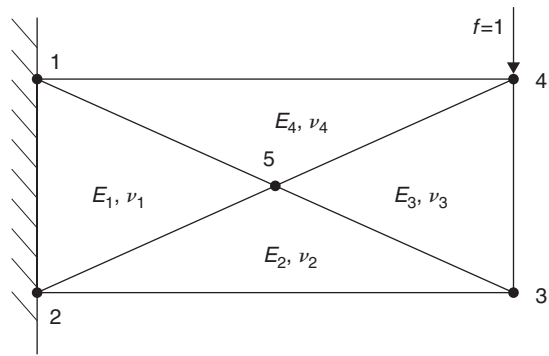

**FIGURE 32.7**  Beam system of finite element inverse analysis.

### 32.3.4  Inverse Initial-Value, Boundary-Value Problem Optimization

Inverse initial-value, boundary-value problem (IVBV) is an optimization paradigm in which GAs have been used successfully (Karr et al., 2000). Some conventional algorithms for solving such search optimizations are the trial-and-error iteration methods that usually start with a guessed solution, and the substitution into the partial differential equations and associated boundary conditions to calculate the errors between predicted values and measured or known values at various locations, then the new guessed or improved solutions are obtained by corrections according to the errors. The aim is to minimize the difference or errors, and the procedure stops once the given precision or tolerance criterion is satisfied. In this way, the inverse problem is actually transformed into an optimization problem.

We now use the heat equation and inverse procedure discussed by Karr et al. (2000) as an example to illustrate the IVBV optimization. On a square plate of unit dimensions, the diffusivity $\kappa(x, y)$ varies with locations $(x, y)$. The heat equation and its boundary conditions can be written as:

$$\frac{\partial u}{\partial t} = \nabla \cdot [\kappa(x,y)\nabla u], \quad 0 < x, y < 1, \quad t > 0,$$
$$u(x, y, 0) = 1, \quad u(x, 0, t) = u(x, 1, t) = u(0, y, t) = u(1, y, t) = 0.$$

The domain is discretized as an $N \times N$ grid, and the measurements of values at $(x_i, y_j, t_n), (i, j = 1, 2, \ldots, N; n = 1, 2, 3)$. The data set consists of the measured value at $N^2$ points at three different times $t_1, t_2, t_3$. The objective is to inverse or estimate the $N^2$ diffusivity values at the $N^2$ distinct locations. The Karr's error metrics are defined as

$$E_u = A \frac{\sum_{i=1}^{N} \sum_{j=1}^{N} \left| u_{i,j}^{\text{measured}} - u_{i,j}^{\text{computed}} \right|}{\sum_{i=1}^{N} \sum_{j=1}^{N} \left| u_{i,j}^{\text{measured}} \right|}, \quad E_\kappa = A \frac{\sum_{i=1}^{N} \sum_{j=1}^{N} \left| \kappa_{i,j}^{\text{known}} - \kappa_{i,j}^{\text{predicted}} \right|}{\sum_{i=1}^{N} \sum_{j=1}^{N} \left| \kappa_{i,j}^{\text{known}} \right|},$$

where $A = 100$ is just a constant.

The floating-point GA proposed by Karr et al. for the inverse IVBV optimization can be summarized as the following procedure: (1) Generate randomly a population containing $N$ solutions to the IVBV problem; the potential solutions are represented in a vector due to the variation of diffusivity $\kappa$ with locations; (2) The error metric is computed for each of the potential solution vectors; (3) $N$ new potential solution vectors are generated by genetic operators such as crossover and mutations in GAs. Selection of solutions depends on the required quality of the solutions with the aim to minimize the error metric and thereby remove solutions with large errors; (4) the iteration continues until the best acceptable solution is found.

On a grid of $16 \times 16$ points with a target matrix of diffusivity, after 40,000 random $\kappa(i, j)$ matrices were generated, the best value of the error metrics $E_u = 4.6050$ and $E_\kappa = 1.50 \times 10^{-2}$. Figure 32.8 shows the error metric $E_\kappa$ associated with the best solution determined by the GA. The small values in the error metric imply that the inverse diffusivity matrix is very close to the true diffusivity matrix.

With some modifications, this type of GA can also be applied to the inverse analysis of other problems such as the inverse parametric estimation of Poisson equation and wave equations. In addition, nonlinear problems can also be studied using GAs.

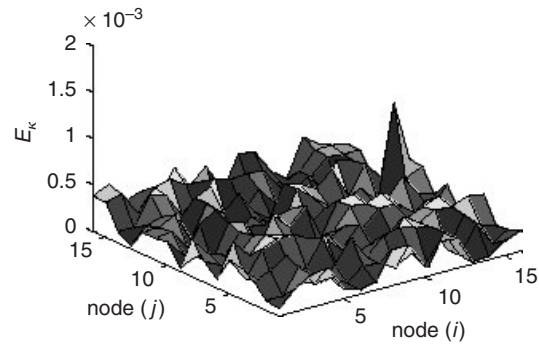

**FIGURE 32.8**　Error metric $E_\kappa$ associated with the best solution obtained by the GA algorithm.

The optimization methods using biology-derived algorithms and their engineering applications have been summarized. We used four examples to show how GAs and PAs can be applied to solve optimization problems in multilevel function optimization, shape design of pressure vessels, finite element inverse analysis of material properties, and the inversion of diffusivity matrix as an IVBV problem. Biology-inspired algorithms have many advantages over traditional optimization methods such as hill-climbing and calculus-based techniques due to parallelism and the ability to locate the best approximate solutions in very large search spaces. Furthermore, more powerful and flexible new generation algorithms can be formulated by combining existing and new evolutionary algorithms with classical optimization methods.

AQ: \* marked References are not cited in text. Please advice.

AQ: Please provide publication details for Heitkotter & Beasley (2000)

AQ: please provide publisher details for Rechenberg (1973)